\documentclass[11pt,titlepage]{article}

\usepackage{rotating, rotfloat}
\usepackage[round]{natbib}
\usepackage{natbib, graphics, rotating, graphicx, hyperref, url, amsmath, amsthm, amsfonts, enumerate, epstopdf, booktabs,verbatim, subfigure, oubraces, array, multirow}
\hypersetup{
    pdfstartview={FitH},    
    colorlinks=true,        
    linkcolor=blue,         
    citecolor=blue,         
    filecolor=blue,         
    urlcolor=blue           
}

     \def\etal{{\itshape et al.}}
\def\beqn{\begin{eqnarray}} \def\eeqn{\end{eqnarray}} 

\def\u{\textbf{\emph{u}}} \def\x{\textbf{\emph{x}}}  \def\z{\textbf{\emph{z}}}  \def\vv{\textbf{\emph{v}}}
\def\y{\textbf{\emph{y}}}

\newcommand{\bec}{\begin{center}}
\newcommand{\enc}{\end{center}}
\newcommand{\bee}{\begin{eqnarray*}}
\newcommand{\ene}{\end{eqnarray*}}
\newcommand{\beq}{\begin{equation}}
\newcommand{\eeq}{\end{equation}}

\begin{document}

\title{\bf The limit of finite sample breakdown point of Tukey's halfspace median for general data}

\author {{Xiaohui Liu$^{a, b}$, Shihua Luo$^{a, b}$, Yijun Zuo$^{c}$}\\ \\
        {\em\footnotesize $^a$ School of Statistics, Jiangxi University of Finance and Economics, Nanchang, Jiangxi 330013, China}\\
         {\em\footnotesize $^b$ Research Center of Applied Statistics, Jiangxi University of Finance and Economics, Nanchang,}\\ {\em\footnotesize Jiangxi 330013, China}\\
         {\em\footnotesize $^c$ Department of Statistics and Probability, Michigan State University, East Lansing, MI, 48823, USA}\\
}

\maketitle

\begin{center}
{\sc Summary}
\end{center}

Under special conditions on data set and underlying distribution, the limit of finite sample breakdown point of
Tukey's halfspace median ($\frac{1} {3}$) has been obtained in literature. In this paper, we establish the result under \emph{weaker assumption} imposed on underlying distribution (halfspace symmetry) and on data set (not necessary in general position).  The representation of Tukey's sample depth regions for data set \emph{not necessary in general position} is also obtained,
as a by-product of our derivation.

\vspace{2mm}

\vspace{2mm}


\setlength{\baselineskip}{1.5\baselineskip}

\vskip 0.1 in
\section{Introduction}
\paragraph{}
\vskip 0.1 in \label{Introduction}

To order multidimensional data, \cite{Tuk1975} introduced the notion of halfspace depth. The  halfspace depth of a point $\x$ in $\mathcal{R}^d$
($d \ge 1$) is defined as
\begin{eqnarray*}
  D(\x, F_n) = \inf_{\u \in \mathcal{S}^{d - 1}} P_n(\u^\top X \leq \u^\top \x),
\end{eqnarray*}
where $\mathcal{S}^{d - 1} = \{\vv\in \mathcal{R}^d: \|\vv\| = 1\}$ with $\|\cdot\|$ being the Euclidean distance, $F_n$ denotes the empirical distribution related to the random sample $\mathcal{X}^n = \{X_1, X_2, \cdots, X_n\}$ from $X \in \mathcal{R}^d$, and $P_n$ is the corresponding empirical probability measure.
\medskip

With this notion, a natural definition of multidimensional median is the point with maximum halfspace depth, which is called
\emph{Tukey's halfspace median} (\emph{HM}). To avoid the nonuniqueness, \emph{HM} ($\hat{\theta}_n$) is defined to be the average of all points lying in the \emph{median region} $\mathcal{M}(\mathcal{X}^n)$, i.e.,
\begin{eqnarray*}
    \hat{\theta}_n := T^*(\mathcal{X}^n) = \textbf{Ave}\left\{\x: \x \in  \mathcal{M}(\mathcal{X}^n) \right\},
\end{eqnarray*}
where $\mathcal{M}(\mathcal{X}^n) = \{\x \in \mathcal{R}^n: D(\x, F_n) = \sup_{\z \in \mathcal{R}^d} D(\z, F_n)\}$, which is the inner-most region among all \emph{$\tau$-trimmed depth regions}:
\begin{eqnarray*}
  \mathcal{D}_\tau(\mathcal{X}^n) = \left\{\x \in \mathcal{R}^d: D(\x, F_n) \ge \tau \right\}, \quad \text{for } \forall \tau \in (0, \lambda^*] \text{ with } \lambda^* = D(\hat{\theta}_n, F_n).
\end{eqnarray*}

When $d = 1$, \emph{HM} reduces to the ordinary univariate median, the latter has the most outstanding property, its best
breakdown robustness. A nature question then is: will \emph{HM} inherit the  best robustness of the univariate median?
\medskip

Answers to this question have been given in the literature, e.g. \cite{DG1992}, \cite{Che1995} and \cite{CT2002} and \cite{AY2002}. The latter two obtained the asymptotic breakdown point ($\frac{1}{3}$) under the maximum bias framework, whereas the former two obtained the limit of finite sample breakdown point (as $n\to \infty$) under the assumption of \emph{absolute continuity} and \emph{central or angular symmetry} of underlying distribution.
\medskip

Among many gauges of robustness of location estimators, finite sample breakdown point is the most prevailing quantitative assessment.
Formally, for a given sample $\mathcal{X}^n$ of size $n$ in $\mathcal{R}^d$, the finite sample addition breakdown point of an location estimator $T$ at $\mathcal{X}^n$ is defined as:
\begin{eqnarray*}
    \varepsilon(T, \mathcal{X}^n) = \min_{1 \leq m \leq n} \left\{\frac{m}{n + m} : \sup_{\mathcal{Y}^m} \left\|T(\mathcal{X}^n \cup \mathcal{Y}^m) -  T(\mathcal{X}^n)\right\| = \infty \right\},
\end{eqnarray*}
where $\mathcal{Y}^m$ denotes a data set of size $m$ with arbitrary values, and $\mathcal{X}^n \cup \mathcal{Y}^m$ the contaminated sample by adjoining $\mathcal{Y}^m$ to $\mathcal{X}^m$.

\medskip
Absolutely continuity guarantees the data set is in \emph{general position} (no more than $d$ sample points lie on a $(d-1)$-dimensional hyperplane \citep{MLB2009}) almost surely.
In practice, the data set $\mathcal{X}^n$ is most likely \emph{not} {in general position}. This is especially true when we are considering the contaminated data set.
\medskip

Unfortunately, most discussions
in literature on finite sample breakdown point is under the assumption of data set in general position. Dropping this unrealistic assumption is very much desirable in the discussion. In this paper we achieve this. Furthermore, we also
relax the angular symmetry \citep{Liu1988, Liu1990} assumption in \cite{Che1995} to a weaker version of symmetry: \emph{halfspace symmetry} \citep{ZS200b}.
$X\in \mathcal{R}^d$ is halfspace symmetrical about $\theta_0$ if $P(X \in \mathcal{H}_{\theta_0}) \ge 1/2$ for any halfspace $\mathcal{H}_{\theta_0}$ containing $\theta_0$. Minimum symmetry is required to guarantee the uniqueness of underlying center $\theta$ in $\mathcal{R}^d$.
\medskip


Without the `in general position' assumption, deriving the limit of finite sample breakdown point of \emph{HM} is quite challenging. We will consider this issue under the combination of halfspace symmetry and a \emph{weak smooth condition} (see Section 2 for details).
Recently,  \cite{LZW2015} have derived the \emph{exact} finite sample breakdown point for fixed $n$. Their result nevertheless depends on the assumption that $\mathcal{X}^n$ is in general position and could \emph{not} be directly utilized under the current setting, because when the underlying $F$ only satisfies the weak smooth condition, the random sample $\mathcal{X}^n$ generated from $F$ may not be in general position in some scenarios. Hence, we have to extend \cite{LZW2015}'s results.
\medskip

Our proofs in this paper heavily depend on the representation of halfspace median region while the existing one in the literature is for the data set in general position. Hence, we have to establish the representation of Tukey's depth region (as the intersection of a finite set of halfspaces) without in general position assumption, which is a byproduct of our proofs.
We only need $\mathcal{X}^n$ to be of affine dimension $d$ which is much \emph{weaker} than the existing ones in \citep{PS2011}. 
\medskip

The rest paper is organized as follows. Section 2 presents a weak smooth condition and shows it is weaker than the absolute continuity and  the interconnection with other notions. Section 3 establishes the representation of Tukey's sample depth regions without in-general-position assumption. Section 4 derives the limiting breakdown point of \emph{HM}. Concluding remarks end the paper.

\vskip 0.1 in
\section{A weak smooth condition}
\paragraph{}
\vskip 0.1 in

In this section, we first present the definition of smooth condition (\emph{SC}), and then investigate its relationship with some other conditions, i.e., absolute continuity and continuous support, commonly assumed in the literature dealing with \emph{HM}. The connection between \emph{SC} and the continuity of the population version of Tukey's depth function $D(\x, F)$ is also investigated.
\medskip

Let $P$ be the probability measure related to $F$. We say a probability distribution $F$ in $\mathcal{R}^d$ of a random vector $X$ is \emph{smooth} at $\x_0 \in \mathcal{R}^d$ if $P(X \in \partial \mathcal{H}) = 0$ for any halfspace $\mathcal{H}$ with $\x_0$ on its boundary $\partial \mathcal{H}$. $F$ is \emph{globally smooth} over $\mathcal{R}^d$ if $F$ is smooth at $\forall \x \in \mathcal{R}^d$.
\medskip

Recall that a distribution $F$ is \emph{absolutely continuous} over $\mathcal{R}^d$ if for $\forall \varepsilon > 0$ there is a positive number $\delta$ such that $P(X \in A) < \varepsilon$ for all Borel sets $A$ of Lebesgue measure less than $\delta$. One can easily show that \emph{absolute continuity} implies \emph{global smoothness}. Nevertheless, the vice versa is \emph{false}. The counterexample can be found in the following.
\medskip

\textbf{Counterexample}. Let $S_1 = \{\x\in \mathcal{R}^d: \|\x\| \leq 1\}$, $S_2 = \{\x\in \mathcal{R}^d: \|\x\| = 2\}$, and $Y = \eta Z_1 + (1 - \eta) Z_2$, where $\eta \sim Bernoulli(0.5)$, and $\eta$, $Z_1$, $Z_2$ are mutually independent. If $Z_1$, $Z_2 \in \mathcal{R}^d$ are uniformly distributed over $S_1$ and $S_2$, respectively, then it is easy to show that the distribution of $Y$ is \emph{not} absolutely continuous, but smooth at $\forall \x \in \mathcal{R}^d$.
\medskip

Furthermore, observe that a distribution $F$ is said to have \emph{contiguous support} if there is no intersection of any two halfspaces with parallel boundaries that has nonempty interior but zero probability and divides the support of $F$ into two parts (see \cite{KZ10}). We can derive that if $F$ has contiguous support it should be globally smooth, but once again the vice versa is \emph{false}. Counterexamples can easily be constructed by following a similar fashion to the above one.
\medskip

Global smoothness is a quite desirable sufficient condition on $F$ if one desires the global continuity of $D(\x, F)$ as shown in the following lemma.
\medskip

\textbf{Lemma 1}. If $F$ is globally smooth, then $D(\x, F)$ is globally continuous in $\x$ over $\mathcal{R}^d$.
\medskip

\textbf{Proof}. When $F$ is globally smooth, we now show that if there $\exists \x_0 \in \mathcal{R}^d$ such that $\lim\limits_{\x \rightarrow \x_0} D(\x, F) \neq D(\x_0, F)$, then it will lead to a contradiction.
\medskip

By noting $\lim\limits_{\x \rightarrow \x_0} D(\x, F) \neq D(\x_0, F)$, we claim that there must exist a sequence $\{\x_k\}_{k = 1}^\infty$ such that $\lim\limits_{k \rightarrow \infty} \x_k = \x_0$ but $\lim\limits_{k \rightarrow \infty} D(\x_k, F) = d_* \neq D(\x_0, F)$. (If $\lim\limits_{k \rightarrow \infty} D(\x_k, F)$ is divergent, by observing $\{D(\x_k, F)\}_{k = 1}^\infty \subset [0, 1]$, we utilize one of its convergent subsequence instead.) For simplicity, hereafter denote $d_k = D(\x_k, F)$ for $k = 0, 1, \cdots $, and assume $d_* < d_0$ if no confusion arises.
\medskip

Observe that $\mathcal{S}^{d - 1}$ is compact. Hence, for each $\x_k$, there $\exists \u_k \in \mathcal{S}^{d - 1}$ satisfying $P(\u_k^\top X \leq \u_k^\top \x_k) = d_k$. Since $\{\u_k\}_{k = 1}^\infty \subset \mathcal{S}^{d - 1}$ is bounded, it should contain a convergent subsequence $\{\u_{k_l}\}_{l = 1}^\infty$ with $\lim\limits_{l \rightarrow \infty} \u_{k_l} = \u_0$. For this $\u_0$ and $\forall \varepsilon_0 \in \left(0,  \frac{d_0 - d_*}{2}\right)$, there $\exists \delta_0 > 0$ such that
\begin{eqnarray}
\label{PBB}
    P(X \in \mathbf{B}(\u_0, \delta_0)) < \varepsilon_0
\end{eqnarray}
following from the global smoothness. Here $\mathbf{B}(\u, c) = \{\z \in \mathcal{R}^d: \u^\top \x_0 - c < \u^\top z \leq \u^\top \x_0\}$ for $\forall \u \in \mathcal{S}^{d-1}$ and $\forall c \in \mathcal{R}^1$.
\medskip

On the other hand, $\lim\limits_{l \rightarrow \infty} P(\u_{k_l}^\top X \leq \u_{k_l}^\top \x_{k_l}) = d_* < d_0 \leq P(\u_{k_l}^\top X \leq \u_{k_l}^\top \x_{0})$. Using this and the convergence of both $\{\x_{k_l}\}_{l = 1}^\infty$ and $\{\u_{k_l}\}_{l = 1}^\infty$, an element derivation leads to that: For $\delta_0$ given above, there $\exists M>0$ such that
\begin{eqnarray*}
  P(X \in \mathbf{B}(\u_{k_l}, \delta_0)) \ge P(\u_{k_l}^\top X \in (\u_{k_l}^\top \x_{k_l}, \u_{k_l}^\top \x_0]) > \frac{d_0 - d_*}{2} > 0,\quad \text{for }\forall k_l > M.
\end{eqnarray*}
This clearly will contradict with \eqref{PBB}, because $\mathbf{B}(\u_0, \delta_0) \cap \mathbf{B}(\u_{k_l}, \delta_0) \rightarrow \mathbf{B}(\u_0, \delta_0)$ as $\u_{k_l} \rightarrow \u_0$ when $k_l \rightarrow \infty$.  \hfill $\Box$
\medskip

Lemma 1 indicates that, when $F$ is globally smooth, $D(\x, F)$ should be globally continuous over $\mathcal{R}^d$, but the vice versa is not clear. Fortunately, an equivalent relationship between the smoothness of $F$ and the continuity of $D(\x, F)$ can be achieved at a special point as stated in the following lemma.
\medskip

\textbf{Lemma 2}. When $F$ is halfspace symmetrical about $\theta_0$, then the following statements are equivalent:
\begin{enumerate}
  \item[(i)] $F$ is smooth at $\theta_0$;

  \item[(ii)] $D(\x, F)$ is continuous at $\theta_0$ with respect to $\x$.
\end{enumerate}

\textbf{Proof}. Similar to Lemma 1, one can show: (i) is true $\Rightarrow$ (ii) is true. In the following, we will show: (i) is false $\Rightarrow$ (ii) is also false.
\medskip

If (i) is false, we claim that there exists a halfspace $\mathcal{H}$ such that $m^0 := P(X \in \partial \mathcal{H}) > 0$. Denote $m^+ = P(X \in \mathcal{H} \setminus \partial \mathcal{H})$ and $m^- = P(X \in \mathcal{H}^c)$. Without confusion, assume that $m^- \leq m^+$ and the normal vector $\u_0$ of $\partial \mathcal{H}$ points into the interior of $\mathcal{H}$. Observe that $D(\x, F) \leq P(\u_0^\top X \leq \u_0^\top \x) \leq \frac{1 - m^0}{2} < 1/2$ for $\forall \x \in \mathcal{H}^c$, i.e., the complementary of $\mathcal{H}$. Hence, for any sequence $\{\x_k\}_{k = 1}^\infty \subset \mathcal{H}^c$ such that $\lim\limits_{k \rightarrow \infty} \x_k = \theta_0$, we have $\limsup\limits_{k \rightarrow \infty} D(\x_k, F) \leq \frac{1 - m^0}{2}$. This in turn implies that $D(\x, F)$ is discontinuous at $\theta_0$ because $D(\theta_0, F) \ge 1/2$ when $F$ is halfspace symmetrical about $\theta_0$.  \hfill $\Box$
\medskip

Summarily, relying on the discussions above, we obtain the following relationship schema when $F$ is halfspace symmetrical about $\theta_0$.
\begin{eqnarray*}
\hline
\left|
  \begin{array}{cl}
  \text{Absolute continuity} &\\
      &\\
      &\\
  \text{Continuous support} &\\
    &
  \end{array}
  \begin{array}{cccc}
  \quad \rotatebox{-5}{\parbox{5mm}{\multirow{3}{*}{$\Longrightarrow$}}}& & &\\
  &F\text{ is globally smooth} &\quad \Rightarrow \quad & F \text{ is smooth at } \theta_0 \\
  \rotatebox{35}{\parbox{5mm}{\multirow{3}{*}{$\Longrightarrow$}}}& & &\\
  &\begin{rotate}{90}$\Leftarrow$\end{rotate}& &\begin{rotate}{90} $\Leftrightarrow$ \end{rotate} \\
  & & &\\
  &D(\x, F) \text{ is globally continuous}&\quad \Rightarrow \quad &\quad D(\x, F) \text{ is continuous at }\theta_0
  \end{array}
\right|
\\
\hline
\end{eqnarray*}
\noindent Since the assumption that $F$ is smooth at $\theta_0$ is quite general, we call it \textbf{\emph{weak smooth condition}} throughout this paper.

\vskip 0.1 in
\section{Representation of Tukey's depth regions}
\paragraph{}
\vskip 0.1 in

To prove the main result, we need to know the representation of $\mathcal{M}(\mathcal{X}^n)$. When $\mathcal{X}^n$ is in general position, this issue has been considered by \cite{PS2011}. Nevertheless, their result can not be directly applied to prove our main theorem, because when the underlying distribution $F$ only satisfies the weak smooth condition, the sample $\mathcal{X}^n$ may \emph{not be in general position}. Hence, we have to solve this problem before proceeding further.
\medskip

For convenience, we introduce the following notations. For $\forall \u \in \mathcal{S}^{d - 1}$ ($d \ge 2$) and $\forall \tau \in (0, 1)$, denote the $(\tau, \u)$-halfspace as
\begin{eqnarray*}
  \mathcal{H}_\tau(\u) = \left\{\x \in \mathcal{R}^d: \u^\top \x \ge q_{\tau}(\u)\right\}
\end{eqnarray*}
with complementary $\mathcal{H}_\tau^c(\u) = \{\x \in \mathcal{R}^d: \u^\top \x < q_{\tau}(\u)\}$ and boundary $\partial \mathcal{H}_\tau(\u) = \{\x \in \mathcal{R}^d: \u^\top \x = q_{\tau}(\u)\}$, where $q_{\tau}(\u) = \inf\{t\in \mathcal{R}^1: F_{\u n}(t) \ge \tau\}$, and $F_{\u n}$ denotes the empirical distribution of $\{\u^\top X_1$, $\u^\top X_2$, $\cdots, \u^\top X_n\}$. 
Obviously, $\u$ points into the interior of $\mathcal{H}_\tau(\u)$, and (see e.g. \cite{KM2008})
\begin{eqnarray}
\label{DTau}
  \mathcal{D}_\tau(\mathcal{X}^n) = \bigcap_{\u\in\mathcal{S}^{d - 1}} \mathcal{H}_\tau(\u).
\end{eqnarray}

\begin{figure}[H]
\centering
    \includegraphics[angle=0,width=4.5in]{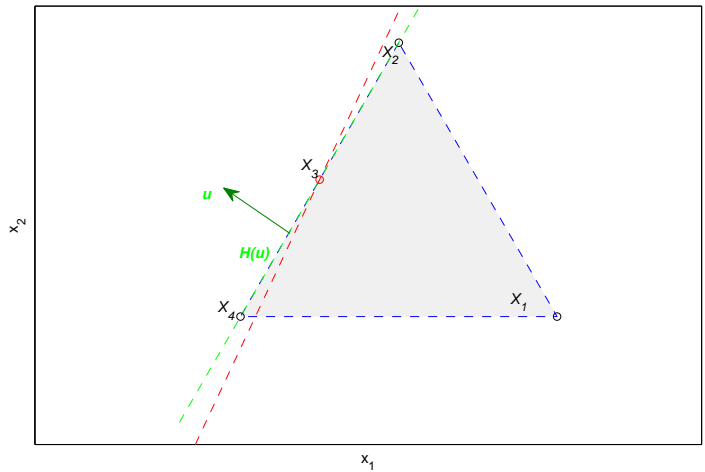}
    \caption{Shown is an example of the $\tau$-\emph{irrotatable} halfspace. Observe that: (a) $\mathcal{H}(\textbf{u})$ cuts away no more than 1 sample point, i.e., $X_1$, and (b) $\partial \mathcal{H}(\u)$ passes through at least $2\ (=d)$ sample points, i.e., $X_2$, $X_3$, $X_4$, and it is possible to make $\mathcal{H}(\u)$ cutting away \emph{more than} $\lceil 4 \times 1/2 \rceil - 1 = 1$ sample points, i.e., $X_1$ and $X_2$, through deviating it by an arbitrary small scale around $X_3$. Hence, $\mathcal{H}(\u)$ is $1/2$-\emph{irrotatable}.}
    \label{fig:tauHalfspace}
\end{figure}

In the following, a halfspace $\mathcal{H}_\tau(\u)$ is said to be \textbf{$\tau$-\emph{irrotatable}} if:
\begin{enumerate}
  \item[(\textbf{a})] $n P_n(X \in \mathcal{H}_\tau^c(\u)) \leq \lceil n \tau \rceil - 1$, i.e., \emph{$\mathcal{H}_\tau(\u)$ cuts away at most} $\lceil n \tau \rceil - 1$ sample points.

  \item[(\textbf{b})]  
        $\partial \mathcal{H}_\tau(\u)$ contains \emph{at least} $d$ sample points, and among them \emph{there exist $d - 1$ points}, which can determine a $(d - 2)$-dimensional hyperplane $\mathbf{V}_{d - 2}$ such that: it is possible to make $\mathcal{H}_\tau(\u)$ cutting away \emph{more than} $\lceil n \tau \rceil - 1$ sample points only through deviating it around $\mathbf{V}_{d - 2}$ by an arbitrary small scale.
\label{PII}
\end{enumerate}
Here $\lceil \cdot \rceil$ denotes the ceiling function, and $\mathbf{V}_{d - 2}$ is a singleton if $d = 2$. To gain more insight, we provide a 2-dimensional example in Figure \ref{fig:tauHalfspace}. In this example, $X_1$, $X_2$, $X_3$ and $X_4$ are clearly not in general position, and $\mathcal{H}(\u)$ is $1/2$-\emph{irrotatable}.

\emph{Remarkably},  if a $\tau_1$-\emph{irrotatable} halfspace cuts away strictly less than $\lceil n \tau_1 \rceil - 1$ sample points, it also should must be $\tau_2$-\emph{irrotatable} for some $\tau_2 < \tau_1$.
\medskip

This $\tau$-\emph{irrotatable} property is quite important for the following lemma, which further plays a key role in the proof of Lemma 4.
\medskip

\textbf{Lemma 3}. Suppose $\mathcal{X}^n = \{X_1, X_2, \cdots, X_n\} \subset \mathcal{R}^d ~ (d \ge 2)$ is of affine dimension $d$. Then for $\forall \tau \in (0, \lambda^*]$, we have
\begin{eqnarray*}
  \mathcal{D}_\tau(\mathcal{X}^n) = \bigcap_{i=1}^{m_\tau} \mathcal{H}_\tau(\mu_i),
\end{eqnarray*}
where $m_\tau$ denotes the number of all $\tau$-\emph{irrotatable} halfspaces $\mathcal{H}_\tau(\mu_i)$.
\medskip

\textbf{Proof}. By \eqref{DTau}, $\mathcal{D}_\tau(\mathcal{X}^n) \subset \bigcap_{i=1}^{m_\tau} \mathcal{H}_\tau(\mu_i)$ holds trivially. Hence, in the sequel we only prove: $\mathcal{D}_\tau(\mathcal{X}^n) \supset \bigcap_{i=1}^{m_\tau} \mathcal{H}_\tau(\mu_i)$.
\medskip

If there $\exists \x_0 \in \bigcap_{i=1}^{m_\tau} \mathcal{H}_\tau(\mu_i)$ such that $\x_0 \notin \mathcal{D}_\tau(\mathcal{X}^n)$, i.e., $D(\x_0, F_n) < \tau$, we now show that this will lead to a contradiction. For simplicity, hereafter denote $\mathcal{V}_n(\tau) = \bigcap_{i=1}^{m_\tau} \mathcal{H}_\tau(\mu_i)$.
\medskip

Since $P_n(\cdot)$ takes values only on $\{0, 1/n, 2/n, \cdots, n/n\}$, there $\exists \u_0 \in \mathcal{S}^{d-1}$ such that
\begin{eqnarray}
\label{u0}
  P_n(\u_0^\top X \leq \u_0^\top \x_0) = D(\x_0, F_n).
\end{eqnarray}
Trivially, when $\mathcal{X}^n$ is of affine dimension $d$, we have: $\mathcal{V}_n(\tau) \subset \textbf{cov}(\mathcal{X}^n)$ for $\forall \tau \in (0, \lambda^*]$, where $\textbf{cov}(\mathcal{X}^n)$ denotes the convex hull of $\mathcal{X}^n$. Hence, for $\x_0 \in \mathcal{V}_n(\tau)$ and $\u_0$ given in \eqref{u0}, there must exist an integer $k_0 \in \{1, 2, \cdots, n\lambda^*\}$ and a permutation $\pi_0 := (i_1, i_2, \cdots, i_n)$ of $(1, 2, \cdots, n)$ such that
\begin{eqnarray}
\label{PermU0}
  \u_0^\top X_{i_1} \leq \u_0^\top X_{i_2} \leq \cdots \leq \u_0^\top X_{i_{k_0}} \leq \u_0^\top \x_0 < \u_0^\top X_{i_{k_0+1}} \leq \cdots \leq \u_0^\top X_{i_{n}}.
\end{eqnarray}
Obviously, $k_0 / n < \tau$ due to $D(\x_0, F_n) < \tau$, and hence $k_0 \leq \lceil n \tau \rceil - 1$.
\medskip

Note that replacing $\u \in \mathcal{S}^{d - 1}$ with $\u \in \mathcal{R}^d \setminus \{0\}$ does no harm to the definition of both $D(\x, F_n)$ and $\mathcal{D}_\tau(\mathcal{X}^n)$ \citep{LZ2014}. Hence, in the sequel we pretend that the constraint on $\u$ is $\u \in \mathcal{R}^d \setminus \{0\}$ instead.
\medskip

Denote $\mathcal{C}(\pi_0) = \{\vv \in \mathcal{R}^d \setminus \{0\}: \vv^\top X_{i_t} \leq \vv^\top X_{i_{k_0 + 1}} \text{ for any } 1\leq t \leq k_0, \text{ and } \vv^\top X_{i_{k_0 + 1}} \leq \vv^\top X_{i_s} \text{ for any } k_0 + 2 \leq s \leq n\}$. Obviously, $\u_0 \in \mathcal{C}({\pi_0})$, and $\mathcal{C}(\pi_0)$ is a convex cone.
\medskip

Let $\mathcal{U} := \{\nu_j\}_{j=1}^{m_v} = \{\z \in \mathcal{R}^d \setminus \{0\}: \|\z\| = 1,~\z \text{ lies in a vertex of } \mathcal{C}({\pi_0})\}$ with $m_v$ being $\mathcal{U}$'s cardinal number. Clearly, $m_v < \infty$ and $\nu_1, \nu_2, \cdots, \nu_{m_\tau}$ are non-coplanar when $\mathcal{X}^n$ is of affine dimension $d$. By the construction of $\mathcal{C}({\pi_0})$, each $\nu \in \mathcal{U}$ determines a halfspace $\mathcal{H}(\nu)$ such that: (\textbf{p1}) $\nu$ is normal to $\partial \mathcal{H}(\nu)$ and points into the interior of $\mathcal{H}(\nu)$, (\textbf{p2}) $\mathcal{H}(\nu)$ cuts away at most $\lceil n\tau \rceil - 1$ sample points, because $X_{i_{k_0 + 1}}, X_{i_{k_0 + 2}}, \cdots, X_{i_{n}} \in \mathcal{H}_\nu$, (\textbf{p3}) $\partial \mathcal{H}(\nu)$ contains \emph{at least} $d$ sample points, which are of affine dimension $d-1$ due to $\nu$ is a vertex of $\mathcal{C}({\pi_0})$. 
\medskip

For $\mathcal{U}$, we claim that: there $\exists \vv_0 \in \mathcal{U}$ satisfying $\vv_0^\top \x_0 < \vv_0^\top X_{i_{k_0 + 1}}$. \emph{If not}, $\nu_j^\top \x_0 \ge \nu_j^\top X_{i_{k_0 + 1}}$ for all $j = 1, 2, \cdots, m_v$. Hence,
\begin{eqnarray*}
    \left(\sum_{j=1}^{m_v} \omega_j \nu_j\right)^\top \x_0 \ge \left(\sum_{j=1}^{m_v} \omega_j \nu_j\right)^\top X_{i_{k_0 + 1}},
\end{eqnarray*}
where $\sum_{j=1}^{m_v} \omega_j = 1$ with $\omega_j \ge 0$ for all $j = 1, 2, \cdots, m_v$. This contradicts with \eqref{PermU0} by noting that $\mathcal{C}(\pi_0)$ is convex and $\u_0 \in \mathcal{C}(\pi_0)$.
\medskip

However, $\vv_0^\top \x_0 < \vv_0^\top X_{i_{k_0 + 1}}$ implies $\x_0 \notin \mathcal{H}(\vv_0)$. We have:
\begin{enumerate}
  \item[\textbf{S1}.] $\mathcal{H}(\vv_0)$ satisfies (\textbf{b}) given in Page \pageref{PII}: By (\textbf{p1})-(\textbf{p3}), $\mathcal{H}(\vv_0)$ is $\tau$-\emph{irrotatable}, contradicting with the definition of $\mathcal{V}_n(\tau)$.

  \item[\textbf{S2}.] $\mathcal{H}(\vv_0)$ does \emph{not} satisfy (\textbf{b}): Among all sample points contained by $\partial \mathcal{H}(\vv_0)$, there must exist $d-1$ points that determine a $(d - 2)$-dimensional hyperplane, around which we can obtain a $\tau$-\emph{irrotatable} halfspace through rotating $\mathcal{H}(\vv_0)$.

(\emph{If not, there will be a contradiction}: By (\textbf{p2}), there $\exists X_{j_1}$, $X_{j_2}$, $\cdots, X_{j_d} \in \partial \mathcal{H}(\vv_0)$, which are of affine dimension $d - 1$. Denote $\mathbf{W}_1, \mathbf{W}_2, \cdots, \mathbf{W}_d$ respectively as ${d \choose d-1}$ hyperplanes that passing through all  $(d-2)$-dimensional facets of the simplex formed by $X_{j_1}, X_{j_2}, \cdots, X_{j_d}$. Then similar to Part (\textbf{II}) of the proof of Theorem 1 in \cite{LLZ2015}, it is easy to check that:
\begin{center}
  for $\forall \y \in \mathcal{R}^d$, $\y$ can \emph{not simultaneously} lie in all $\mathbf{W}_1, \mathbf{W}_2, \cdots, \mathbf{W}_d$.
\end{center}
Without confusion, assume $\y \notin \mathbf{W}_1$ and $X_{j_1} \in \mathbf{W}_1$. Observe that no $\tau$-\emph{irrotatable} halfspace is available through rotating $\mathcal{H}(\vv_0)$ around $\mathbf{W}_1$. Hence, for $\forall \delta > 0$,
    \begin{eqnarray}
       \label{IRhalfspace}
        \max\{nP_n(X \in \mathcal{H}_{\delta+}^c), ~ nP_n(X \in \mathcal{H}_{\delta-}^c) \} < \lceil n \tau \rceil - 1,
    \end{eqnarray}
       where $\mathcal{H}_{\delta+} = \{\z \in \mathcal{R}^d: \u_{+}^\top \z \ge \u_{+}^\top X_{j_1}\}$, and $\mathcal{H}_{\delta-} = \{\z \in \mathcal{R}^d: \u_{-}^\top \z \ge \u_{-}^\top X_{j_1}\}$ with $\u_+ = \vv_0 + \delta \u_*$ and $\u_- = \vv_0 - \delta \u_*$, where $\u_* \in \mathcal{S}^{d - 1}$ is orthogonal to both $\vv_0$ and $\mathbf{W}_1$. Since either $\y \in \mathcal{H}_{\delta+}^c$ or $\y \in \mathcal{H}_{\delta-}^c$ for $\forall \delta > 0$, we obtain $D(\y, F_n) < (\lceil n \tau \rceil - 1)/n \leq \tau$. This is \emph{impossible} because $\mathcal{D}_n(\tau)$ is nonempty for $\forall \tau \in (0, \lambda^*]$.)

Furthermore, it is easy to show that: if there is a $\tau$-\emph{irrotatable} halfspace, say $\mathcal{H}_1$, obtained through rotating $\mathcal{H}(\vv_0)$ around one $(d - 2)$-dimensional hyperplane clockwise (without confision), then there would be an another $\tau$-\emph{irrotatable} halfspace, say $\mathcal{H}_2$, by rotating $\mathcal{H}(\vv_0)$ anti-clockwise. By noting $\mathcal{H}^c(\vv_0) \subset \mathcal{H}_1^c \cup \mathcal{H}_2^c$, we can obtain either $\x_0 \in \mathcal{H}_1^c$ or $\x_0 \in \mathcal{H}_2^c$, which contradicts with the definition of $\mathcal{V}_n(\tau)$.
\end{enumerate}
Hence, there is \emph{no} such $\x_0$ that $\x_0 \in \mathcal{V}_n(\tau)$, but $\x_0 \notin \mathcal{D}_\tau(\mathcal{X}^n)$.
\medskip

This completes the proof of this lemma.  \hfill $\Box$
\medskip

\noindent \textbf{Remark 1}. It may have long been known in the statistical community that Tukey's sample depth regions may be polyhedral and have a finite number of facets. The detailed character of each facet of these regions is \emph{unknown}, nevertheless. When $\mathcal{X}^n$ is in general position, \cite{PS2011} have shown that each hyperplane passing through a facet of $\mathcal{D}_\tau(\mathcal{X}^n)$, for $\forall \tau \in (0, \lambda^*]$, contains \emph{exactly} $d$ and cuts away exactly $\lceil n \tau \rceil - 1$ sample points; see Lemma 4.1 in Page 201 of \cite{PS2011} for details. Lemma 3 generalizes their result by removing the `in general position' assumption, and indicates that such hyperplanes contain \emph{at least} $d$ and cuts away \emph{no more than} $\lceil n \tau \rceil - 1$ sample points.

\begin{figure}[H]
\begin{center}
	\subfigure[Halfspace $\mathcal{H}_{1/2}(u_1)$, which is $1/2$-\emph{irrotatable} because $4P_n(X\in \mathcal{H}_\tau^c(u_1)) \leq \lceil 4\tau \rceil -1$ but $4P_n(X \notin \mathcal{H}_\tau^*(u_1)) = 2 > \lceil 4\tau \rceil -1$ for $\tau = 1/2$.]{
	\includegraphics[angle=0,width=4.2in]{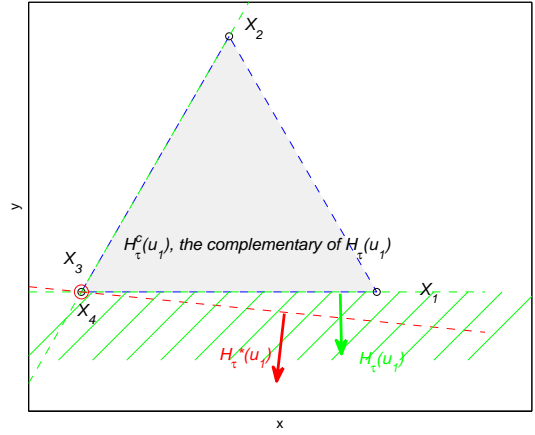}
	\label{fig:Case1}
	}
	\subfigure[Halfspace $\mathcal{H}_{1/2}(u_2)$, which is similarly $1/2$-\emph{irrotatable} because $4P_n(X\in \mathcal{H}_\tau^c(u_2)) \leq \lceil 4\tau \rceil -1$ but $4P_n(X \notin \mathcal{H}_\tau^*(u_2)) = 2 > \lceil 4\tau \rceil -1$ for $\tau = 1/2$]{
	\includegraphics[angle=0,width=4.2in]{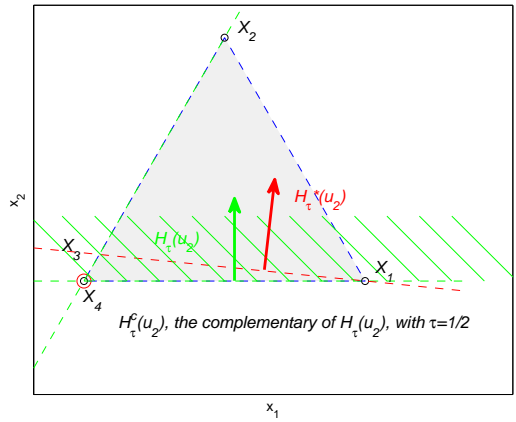}
	\label{fig:Case2}
	}
	\subfigure[The intersection of lines $L_1$ and $L_2$.]{
	\includegraphics[angle=0,width=4.2in]{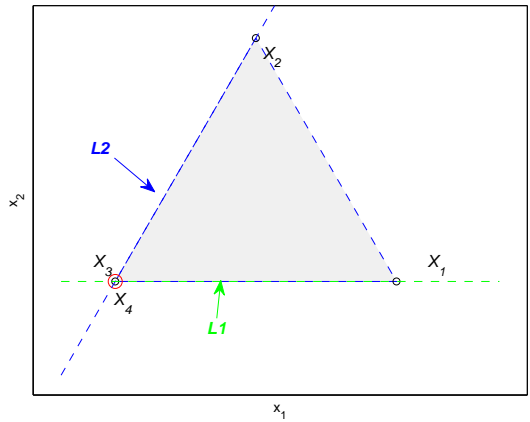}
	\label{fig:Intersection}
	}
\caption{Shown are examples of the $\tau$-\emph{irrotatable} halfspaces and related \emph{HM}.}
\label{fig:DepthRegion}
\end{center}
\end{figure}

To facilitate the understanding, we provide an illustrative example in Figure \ref{fig:DepthRegion}. In this example, there are $n = 4$ observations, i.e., $X_1, X_2, X_3, X_4$, where $X_3$ and $X_4$ take the same value. Clearly, they are \emph{not} in general position and of affine dimension 2. 
Figures \ref{fig:Case1}-\ref{fig:Case2} indicate that $\{X_1, X_3, X_4\}$ determines two $1/2$-\emph{irrotatable} halfspaces, i.e., $\mathcal{H}_{1/2}(\u_1)$ and $\mathcal{H}_{1/2}(\u_2)$, satisfying $\mathcal{H}_{1/2}(\u_1) \cap \mathcal{H}_{1/2}(\u_2) = L_1$. Similarly, the intersection of the halfspaces determined by $\{X_2, X_3, X_4\}$ is $L_2$. Hence, the median region is $\{\x: \x = X_3\}$. \emph{From Figure \ref{fig:Case2} we can see that $\partial \mathcal{H}_{1/2}(\u_2)$ contains $3\ (\neq 2)$ and $\mathcal{H}_{1/2}(\u_2)$ cuts away $0\ (\neq 1)$ sample points, which obviously is \textbf{not} in agreement with the results of \cite{PS2011}.}

\vskip 0.1 in
\section{The limiting breakdown point of \emph{HM}}
\paragraph{}
\vskip 0.1 in

In this section, we will derive the limit of the finite sample breakdown point of \emph{HM} when the underlying distribution satisfies only the weak smooth condition (such a limit is also called asymptotic breakdown point in the literature, the latter notion is based on the maximum bias notion though, see Hampel(1968)). Since \emph{HM} reduces to the ordinary univariate median for $d=1$, whose breakdown point robustness has been well studied, \emph{we focus only on the scenario of $d \ge 2$} in the sequel.

The key idea is to obtain simultaneously a lower and an upper bound of $\varepsilon(T^*, \mathcal{X}^n)$ for fixed $n$, and then prove that they tend to the same value as $n \rightarrow \infty$. When $\mathcal{X}^n$ is of affine dimension $d$, it is easy to obtain a lower bound, i.e., $\frac{\lambda^*}{1 + \lambda^*}$, for $\varepsilon(T^*, \mathcal{X}^n)$ by using a similar strategy to \cite{DG1992} though. \emph{Finding a proper upper bound is not trivial}, nevertheless.
\medskip

To this end, we establish the following lemma, which provides a sharp upper bound with its limit coinciding with that of $\frac{\lambda^*}{1 + \lambda^*}$ asymptotically. For simplicity, denoting by $\mathbb{A}_\u$ an arbitrary $d \times (d - 1)$ matrix of unit vectors such that $(\u\,\vdots\,\mathbb{A}_{\u})$ constitutes an orthonormal basis of $\mathcal{R}^d$, we define the $\mathbb{A}_\u$-projections of $\mathcal{X}^n$ as $\mathbf{X}_\u^n = \{\mathbb{A}_\u^\top X_1, \mathbb{A}_\u^\top X_2, \cdots, \mathbb{A}_\u^\top X_n\}$ for $\forall \u \in \mathcal{S}^{d - 1}$. Correspondingly, let $\hat{\theta}_n^\u = T^*(\mathbf{X}_\u^n)$, $\lambda_\u^* = D(\hat{\theta}_n^\u, F_{\u n})$, and $F_{\u n}$ to be the empirical distribution related to $\mathbf{X}_\u^n$.
\medskip

\textbf{Lemma 4}. For a given data set $\mathcal{X}^n$ of affine dimension $d$, the finite sample breakdown point of Tukey's halfspace median satisfies 
\begin{eqnarray*}
  \varepsilon(T^*, \mathcal{X}^n) \leq \frac{\inf_{{\u\in \mathcal{S}^{d - 1}}} \lambda_\u^*}{1 + \inf_{{\u\in \mathcal{S}^{d - 1}}} \lambda_\u^*}.
\end{eqnarray*}

\begin{figure}[H]
\centering
    \includegraphics[angle=0,width=3.0in]{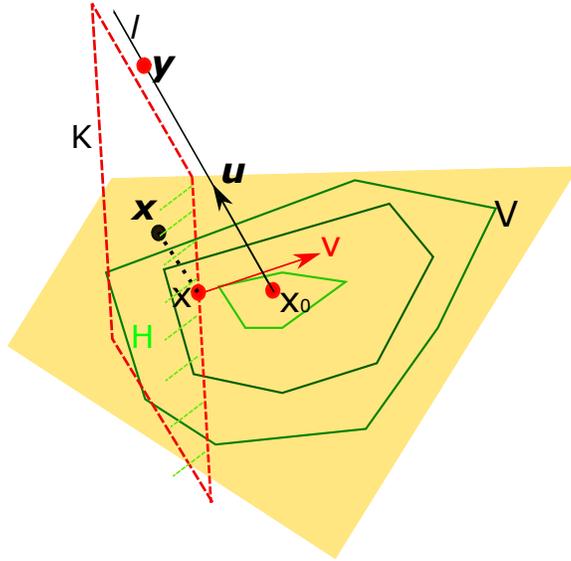}
    \caption{Shown is a 3-dimensional illustration. Once $\y$'s are putted on $\ell$, all of their projections onto $\mathbf{V}$ are $\textbf{x}_0$. Hence, for any $\x \notin \ell$, its depth with respect to $\mathcal{X}^n\cup \mathcal{Y}^m$ would be no more than that of $\mathbf{x}$ with respect to the projections of $\mathcal{X}^n\cup \mathcal{Y}^m$, because all projections of the sample points contained by $\textbf{K}$ would lie in $\mathbf{H}$. Here $\mathbf{H}$ denotes the optimal $(d-1)$-dimensional optimal halfspace of $\textbf{x}$, and $\mathbf{K}$ the $d$-dimensional halfspace whose projection is $\mathbf{H}$.}
    \label{fig:Lemma5}
\end{figure}

Since the whole proof of this lemma is very long, we present it in two parts. For $\forall \u \in \mathcal{S}^{d-1}$, in Part (\textbf{I}), we first project $\mathcal{X}^n$ onto a $(d-1)$-dimensional space $\textbf{V}_{d-1}^\u$ that is orthogonal to $\u$, and then show that there $\exists \textbf{x}_0 \in \textbf{V}_{d-1}^\u$, which can lie in the inner of the complementary of a $(d-1)$-dimensional optimal halfspace of $\forall \textbf{x} \in \textbf{V}_{d-1}^\u \setminus \{\textbf{x}_0\}$. Here by optimal halfspace of $\textbf{x}$ we mean the halfspace realizing the depth at $\textbf{x}$ with respect to $\mathbf{X}_\u^n$. Denote the line passing through $\mathbf{x}_0$ and parallel to $\u$ as $\ell_\u$. In Part (\textbf{II}), we will show that by putting $n \lambda_\u^*$  repetitions of $\y_0$ at any position on $\ell_\u$ but outside the convex hull of $\mathcal{X}^n$, i.e., $\ell_\u \setminus \textbf{cov}(\mathcal{X}^n)$, it is possible to obtain that $\sup_{\x \in \textbf{cov}(\mathcal{X})} D(\x, \mathcal{X}^n \cup \mathcal{Y}^m) \leq n \lambda_\u^*$. Hence, $\inf_{\u \in \mathcal{S}^{d-1}} n \lambda_\u^*$  repetitions of $\y_0$ suffice for breaking down $T^*(\mathcal{X}^n \cup \mathcal{Y}^m)$. See Figure \ref{fig:Lemma5} for a 3-dimensional illustration.
\medskip

\textbf{Proof of Lemma 4}. Trivially, it is easy to check that, for $\forall \u \in \mathcal{S}^{d - 1}$, $\mathbf{X}_\u^n$ is of affine dimension $d-1$ if $\mathcal{X}^n$ is of affine dimension $d$.
\medskip

(\textbf{I}). In this part, we only prove that: When the affine dimension of $\mathcal{M}(\mathbf{X}_\u^n)$ is nonzero for $d > 2$, there $\exists \mathbf{x}_0 \in \mathcal{M}(\mathbf{X}_u^n)$ such that $\mathcal{U}_{\mathbf{x}} \cap \mathcal{H}_{\mathbf{x}, \mathbf{x}_0} \neq \emptyset$ for $\forall \mathbf{x} \in \mathcal{M}(\mathbf{X}_\u^n) \setminus \{\mathbf{x}_0\}$, where $\mathcal{U}_{\mathbf{x}} = \{\mathbf{v} \in \mathcal{S}^{d - 2}: P_n(\mathbf{v}^\top (\mathbb{A}_\u^\top X) \leq \mathbf{v}^\top \mathbf{x}) = D(\mathbf{x}, F_{\u n})\}$, and $\mathcal{H}_{\mathbf{x}, \mathbf{x}_0} = \{\mathbf{v} \in \mathcal{S}^{d - 2}: \mathbf{v}^\top \mathbf{x} < \mathbf{v}^\top \mathbf{x}_0\}$. That is, $\mathbf{x}_0$ lies in the inner of the complementary of a $(d-1)$-dimensional optimal halfspace of $\forall \textbf{x} \in \mathcal{M}(\mathbf{X}_\u^n) \setminus \{\textbf{x}_0\}$. The rest proof follows a similar fashion to Lemmas 2-3 of \cite{LZW2015}.
\medskip

By Lemma 3, $\mathcal{M}(\mathbf{X}_\u^n)$ is polyhedral. Similar to Theorem 2 of \cite{LLZ2015}, we can obtain that, if there is a sample point $X_i$ such that $\mathbb{A}_\u^\top X_i \in \mathcal{M}(\mathbf{X}_\u^n)$, then $\mathbb{A}_\u^\top X_i$ should be a vertex of $\mathcal{M}(\mathbf{X}_\u^n)$ based on the representation of $\mathcal{M}(\mathbf{X}_\u^n)$ obtained in Lemma 3. Let $\mathbf{V}_\u$ be the set of vertexes of $\mathcal{M}(\mathbf{X}_\u^n)$ such that, for $\forall \textbf{y} \in \mathbf{V}_\u$, there is an optimal halfspace $\mathbf{H}_\textbf{y}$ of $\textbf{y}$ satisfying $\mathbf{H}_\textbf{y}  \cap \mathcal{M}(\mathbf{X}_\u^n) = \{\textbf{y}\}$. Trivially, $\mathbb{A}_\u^\top X_i \in \mathbf{V}_\u$ if $\mathbb{A}_\u^\top X_i \in \mathcal{M}(\mathbf{X}_\u^n)$.
\medskip

If there is point in $\mathbf{V}_\u$ that can sever as $\mathbf{x}_0$, then this statement holds already. Otherwise, \emph{find a candidate point $\textbf{z}_0$ by using the following iterative procedure and then show that $\textbf{z}_0$ can be used as} $\textbf{x}_0$. For simplicity, hereafter denote $\mathcal{A}_{\mathbf{z}} = \{\mathbf{x} \in \mathcal{R}^{d - 1}: \mathcal{U}_{\mathbf{x}} \cap \mathcal{H}_{\mathbf{x}, \mathbf{z}} \neq \emptyset\}$ and $\mathcal{B}_{\mathbf{z}} = \{\mathbf{x} \in \mathcal{R}^{d - 1}: \mathcal{U}_{\mathbf{x}} \cap \mathcal{H}_{\mathbf{x}, \mathbf{z}} = \emptyset\}$ for $\forall \mathbf{z} \in \mathcal{M}(\mathbf{X}_u^n)$. Obviously, $\mathcal{A}_{\mathbf{z}} \cup \mathcal{B}_{\mathbf{z}} = \mathcal{R}^{d - 1}$, $\mathcal{A}_{\mathbf{z}} \cap \mathcal{B}_{\mathbf{z}} = \emptyset$, $\mathbf{z} \in \mathcal{B}_{\mathbf{z}}$, and $\mathcal{B}_{\mathbf{z}} \subset \mathcal{M}(\mathbf{X}_\u^n)$.
\medskip

Let $\mathbf{z}_1 = T^*(\mathbf{X}_\u^n)$. Clearly, $\mathbf{V}_\u \cap \mathcal{B}_{\mathbf{z}_1}  = \emptyset$. (In fact, $\mathbf{V}_\u \cap \mathcal{B}_{\mathbf{z}}  = \emptyset$ for any $\mathbf{z} \in \mathcal{M}(\mathbf{X}_\u^n) \setminus \mathbf{V}_\u$.) If $\mathcal{B}_{\mathbf{z}_1} = \{\textbf{z}_1\}$, let $\mathbf{x}_0 = \mathbf{z}_0$ and this statement is already true. Otherwise, similar to Lemma 2 of \cite{LZW2015}, for $\forall \mathbf{x} \in \mathcal{B}_{\mathbf{z}_1} \setminus \{\mathbf{z}_1\}$, we obtain: (\textbf{o1}) $\mathbf{u}^\top \mathbf{x} \ge \mathbf{u}^\top \mathbf{z}_1$ for $\forall \mathbf{u} \in \mathcal{U}_{\mathbf{z}_1}$, (\textbf{o2}) $\mathcal{U}_{\mathbf{x}} \subset \mathcal{U}_{\mathbf{z}_1}$, and (\textbf{o3}) $\mathcal{B}_{\mathbf{x}} \subset \mathcal{B}_{\mathbf{z}_1} \setminus \{\mathbf{z}_1\}$.
\medskip

Denote
\begin{eqnarray*}
  g(\mathbf{z}_1) = \sup_{\mathbf{v} \in \mathcal{U}_{\mathbf{z}_1}, \mathbf{x} \in \mathcal{B}_{\mathbf{z}_1 \setminus \{\mathbf{z}_1\}}} \mathbf{v}^\top (\mathbf{x} - \mathbf{z}_1).
\end{eqnarray*}
Clearly, $g(\mathbf{z}_1) > 0$ by (\textbf{o1})-(\textbf{o3}). Along the same line of \cite{LZW2015}, we can find a series $\{\mathbf{z}_i\}_{i = 1}^{\infty} \subset \mathcal{M}(\mathbf{X}_u^n)$, if there is \emph{no} $m > 1$ such that $\mathcal{B}_{\mathbf{z}_m} = \{\mathbf{z}_m\}$, satisfying that: $\{\mathbf{z}_i\}_{i = 1}^{\infty}$ contains a convergent subsequence $\{\mathbf{z}_{i_k}\}_{k = 1}^{\infty}$ with $\lim\limits_{k \rightarrow \infty} \mathbf{z}_{i_k} = \mathbf{z}_0$ and $\lim\limits_{k \rightarrow \infty} g(\mathbf{z}_{i_k - 1}) = 0$. Trivially, $\mathbf{z}_0 \in \mathcal{M}(\mathbf{X}_\u^n) \setminus \mathbf{V}_\u$. (If not, it is easy to obtain a contradiction.)
\medskip

\emph{Now we proceed to prove $\mathcal{B}_{\mathbf{z}_0} = \{\mathbf{z}_0\}$}. \emph{First}, we show
\begin{center}
  (\textbf{F1}): \quad $\mathbf{z}_0 \in \mathcal{B}_{\mathbf{z}_{j - 1}} \setminus \{\mathbf{z}_{j - 1}\}$ for $\forall j \in \{i_k\}_{k = 1}^\infty$.
\end{center}
If not, there must $\exists \tilde{\mathbf{u}} \in \mathcal{U}_{\mathbf{z}_0}$ satisfying $\tilde{\mathbf{u}}^\top \mathbf{z}_0 < \tilde{\mathbf{u}}^\top \mathbf{z}_{j - 1}$. For this $\tilde{\mathbf{u}} \in \mathcal{U}_{\mathbf{z}_0}$, let $(i_1', i_2', \cdots, i_n')$ be the permutation of $(1, 2, \cdots, n)$ such that: (a) $\tilde{\mathbf{u}}^\top (\mathbb{A}_\u^\top X_{i_s'}) \leq \tilde{\mathbf{u}}^\top \mathbf{z}_0$ for $1 \leq s \leq k^*$, and (b) $\tilde{\mathbf{u}}^\top (\mathbb{A}_\u^\top X_{i_t'}) > \tilde{\mathbf{u}}^\top \mathbf{z}_0$ for $k^* + 1 \leq t \leq n$, where $k^* = n \lambda_\u^*$. Denote
\begin{eqnarray*}
  \varepsilon_0 = \frac{1}{2} \min\left\{\min_{k^* + 1 \leq t \leq n}\tilde{\mathbf{u}}^\top ((\mathbb{A}_\u^\top X_{i_{t}'}) - \mathbf{z}_0), ~ \tilde{\mathbf{u}}^\top (\mathbf{z}_{j - 1} - \mathbf{z}_0)\right\}.
\end{eqnarray*}
Since $\{\mathbf{z}_{i_k}\}_{k = 1}^{\infty}$ is convergent, there must $\exists j^* \in \{i_k\}_{k=1}^\infty$ with $j^* > j$ such that $\|\mathbf{z}_{j^*} - \mathbf{z}_0\| < \varepsilon_0$. This, together with $|\tilde{\mathbf{u}}^\top (\mathbf{z}_{j^*} - \mathbf{z}_0)| \leq \|\mathbf{z}_{j^*} - \mathbf{z}_0\|$, leads to $\tilde{\mathbf{u}}^\top \mathbf{z}_{j^*} < \tilde{\mathbf{u}}^\top (\mathbb{A}_u^\top X_{i_t'})$ for $k^* + 1 \leq t \leq n$, which further implies $P_n(\tilde{\mathbf{u}}^\top (\mathbb{A}_u^\top X) \leq \tilde{\mathbf{u}}^\top \mathbf{z}_{j^*}) \leq \lambda_\u^*$. Next, by noting $\lambda_\u^* = D(\mathbf{z}_{j^*}, F_n^u) \leq P_n(\tilde{\mathbf{u}}^\top (\mathbb{A}_u^\top X) \leq \tilde{\mathbf{u}}^\top \mathbf{z}_{j^*})$, we obtain $\tilde{\mathbf{u}} \in \mathcal{U}_{\mathbf{z}_{j^*}} \subset \mathcal{U}_{\mathbf{z}_{j - 1}}$. On the other hand, for $\varepsilon_0$, a similar derivation leads to $\tilde{\mathbf{u}}^\top \mathbf{z}_{j^*} < \tilde{\mathbf{u}}^\top \mathbf{z}_{j - 1}$. This contradicts with $\mathbf{z}_{j^*} \in \mathcal{B}_{\mathbf{z}_{j - 1}}$ when $j^* > j$ by (\textbf{o1})-(\textbf{o2}). \emph{Then}, based on $\lim\limits_{k \rightarrow \infty} g(\mathbf{z}_{i_k - 1}) = 0$ and (\textbf{F1}), we can obtain $\mathcal{B}_{\mathbf{z}_0} \setminus \{\mathbf{z}_0\} = \emptyset$ similar to Lemma 3 of \cite{LZW2015}. Hence, we may let $\mathbf{x}_0 = \mathbf{z}_0$.
\medskip

(\textbf{II}). By denoting $\ell_\u = \{\z \in \mathcal{R}^d: \z = \mathbb{A}_{\u} \mathbf{x}_0 + \gamma \u, ~ \forall \gamma \in \mathcal{R}^1 \}$ and using a similar method to the first proof part of Theorem 1 in \cite{LZW2015}, we can obtain that, for an any given $\y_0 \in \textbf{cov}(\mathcal{X}^n) \setminus \ell_\u$, it holds $\sup_{\x \in \textbf{cov}(\mathcal{X}^n)} D(\x, F_{n+m}) \leq \frac{n \lambda_\u^*}{n + m}$, where $F_{n+m}$ denotes the empirical distribution related to $\mathcal{X}^n\cup \mathcal{Y}^m$, and $\mathcal{Y}^m$ contains $m$ repetitions of $\y_0$.
\medskip

Note that $\u$ is any given, and $D(\y_0, F_{n+m}) = \frac{m}{n + m} \ge \frac{n \lambda_\u^*}{n + m}$ when $m \leq n \lambda_\u^*$. Hence
\begin{eqnarray*}
  \varepsilon(T^*, \mathcal{X}^n) \leq \frac{\inf_{{\u\in \mathcal{S}^{d - 1}}} n \lambda_\u^*}{n + \inf_{{\u\in \mathcal{S}^{d - 1}}} n \lambda_\u^*} = \frac{\inf_{{\u\in \mathcal{S}^{d - 1}}} \lambda_\u^*}{1 + \inf_{{\u\in \mathcal{S}^{d - 1}}} \lambda_\u^*}.
\end{eqnarray*}
This completes the proof.\hfill $\Box$
\medskip

Observe that the upper bound given in Lemma 4 involves the $\mathbb{A}_\u$-projections. A nature problem arises: whether the $\mathbb{A}_\u$-projection of $X$ is still halfspace symmetrically distributed? The following lemma provides a positive answer to this question.
\medskip

\textbf{Lemma 5}. Suppose $X$ is halfspace symmetrical about $\theta_0 \in \mathcal{R}^d$ ($d \ge 2$). Then for $\forall \u \in \mathcal{S}^{d-1}$, $\mathbb{A}_\u^\top X$ is halfspace symmetrical about $\mathbb{A}_\u^\top \theta_0 \in \mathcal{R}^{d-1}$.
\medskip

\textbf{Proof}. For $\forall \mathbf{v} \in \mathcal{S}^{d-2}$, the fact $(\mathbb{A}_\u \mathbf{v})^\top (\mathbb{A}_\u \mathbf{v}) = \mathbf{v}^\top (\mathbb{A}_\u^\top \mathbb{A}_\u) \mathbf{v} = 1$ implies $\mathbb{A}_\u \mathbf{v} \in \mathcal{S}^{d-1}$. Note that
\begin{eqnarray*}
    P\left(\mathbf{v}^\top (\mathbb{A}_\u^\top X) \ge \mathbf{v}^\top (\mathbb{A}_\u^\top \theta_0)\right) = P\left((\mathbb{A}_\u \mathbf{v})^\top X \ge (\mathbb{A}_\u \mathbf{v})^\top \theta_0\right) \ge \frac{1}{2}.
\end{eqnarray*}
This completes the proof of this lemma.
\medskip

Lemma 5 in fact obtains the population version, i.e., $D(\mathbb{A}_\u^\top \theta_0, F_{\u})$, of $D(\hat{\theta}_n^\u, F_{\u n})$ for $\forall \u \in \mathcal{S}^{d - 1}$, where $F_{\u}$ denotes the distribution of $\mathbb{A}_\u^\top X$.
\medskip

we now are in the position to prove the following theorem.
\medskip

\textbf{Theorem 1}. Suppose that (C1) $\{X_1, X_2, \cdots, X_n\} \stackrel{\text{i.i.d.}}{\sim} F$ is of affine dimension $d$, (C2) $F$ is halfspace symmetric about point $\theta_0$, and (C3) $F$ is smooth at point $\theta_0$. Then we have
$
  \varepsilon(T^*, \mathcal{X}^n) \stackrel{\text{a.s.}}{\longrightarrow} \frac{1}{3}, \quad \text{as } n \rightarrow +\infty,
$
where $\stackrel{\text{a.s.}}{\longrightarrow}$ denotes the ``almost sure convergence".
\medskip

\textbf{Proof}. Observe that
\begin{eqnarray*}
  |D(\hat{\theta}_n, F_n) - D(\theta_0, F)| \leq \underbrace{\sup_{\x \in \mathcal{R}^d} |D(\x, F_n) - D(\x, F)|}_\textbf{E1} + \underbrace{|D(\hat{\theta}_n, F) - D(\theta_0, F)|}_\textbf{E2}.
\end{eqnarray*}

Under Condition (C1), a direct use of Remark 2.5 in Page 1465 of \cite{Zuo2003} leads to that
\begin{eqnarray}
\label{supE1}
  \sup_{\x \in \mathcal{R}^d} \left|D(\x, F_n) - D(\x, F)\right| \stackrel{\text{a.s.}}{\longrightarrow} 0, \quad \text{as } n \rightarrow +\infty,
\end{eqnarray}
holds with no restriction on $F$. Hence, $\textbf{E1} \stackrel{\text{a.s.}}{\longrightarrow} 0$.
\medskip

For \textbf{E2}, from Lemma 2, we have that $D(\x, F)$ is continuous at $\theta_0$ under Condition (C3). On the other hand, since $D(\theta_0, F) \ge 1/2 > 0$ under Condition (C2), an application of Lemma A.3 of \cite{Zuo2003} leads to $\hat{\theta}_n \stackrel{\text{a.s.}}{\longrightarrow} \theta_0$. These two facts together imply $\textbf{E2} \stackrel{\text{a.s.}}{\longrightarrow} 0$.
\medskip

Based on $\textbf{E1} \stackrel{\text{a.s.}}{\longrightarrow} 0$ and $\textbf{E2} \stackrel{\text{a.s.}}{\longrightarrow} 0$, we in fact obtain
\begin{eqnarray}
\label{leftPart}
  D(\hat{\theta}_n, F_n) \stackrel{\text{a.s.}}{\longrightarrow} D(\theta_0, F), \quad \forall \u \in \mathcal{S}^{d - 1}.
\end{eqnarray}
Relying on this and the lower bound $\frac{\lambda^*}{1 + \lambda^*}$, it is easy to show that
\begin{eqnarray*}
  \varepsilon(T^*, \mathcal{X}^n) \ge \frac{D(\theta_0, F)}{1 + D(\theta_0, F)},\quad \text{almost surely.}
\end{eqnarray*}

By Lemma 5, $F_{\u}$ is also halfspace symmetrical and smooth at $\mathbb{A}_\u^\top \theta_0$ for $\forall \u \in \mathcal{S}^{d - 1}$. Hence, a similar proof to \eqref{leftPart} leads to
\begin{eqnarray*}
    D(\hat{\theta}_n^\u, F_{\u n}) \stackrel{\text{a.s.}}{\longrightarrow} D(\mathbb{A}_\u^\top \theta_0, F_{\u}), \quad \text{as } n \rightarrow \infty.
\end{eqnarray*}
This, together with Lemma 4, and the theory of empirical processes \citep{Pol1984}, leads to
\begin{eqnarray*}
  \varepsilon(T^*, \mathcal{X}^n) \leq \frac{\inf_{{\u\in \mathcal{S}^{d - 1}}} D(\mathbb{A}_\u^\top \theta_0, F_{\u})}{1 + \inf_{{\u\in \mathcal{S}^{d - 1}}} D(\mathbb{A}_\u^\top \theta_0, F_{\u})},\quad \text{almost surely.}
\end{eqnarray*}

Next, by
\begin{eqnarray*}
  D(\mathbb{A}_\u^\top \theta_0, F_\u) &=& \inf_{\mathbf{v} \in \mathcal{S}^{d-2}} P\left(\mathbf{v}^\top (\mathbb{A}_\u^\top X) \leq \mathbf{v}^\top (\mathbb{A}_\u^\top \theta_0)\right)\\
  &=& \inf_{\bar{\u} \in \mathcal{S}^{d-1},~\bar{\u} \bot \u} P(\bar{\u}^\top X \leq \bar{\u}^\top \theta_0),
\end{eqnarray*}
where $\bar{\u} = \mathbb{A}_\u \mathbf{v}$, we obtain
\begin{eqnarray*}
  \inf_{\u\in \mathcal{S}^{d - 1}} D(\mathbb{A}_\u^\top \theta_0, F_{\u}) &=& \inf_{\u\in \mathcal{S}^{d - 1}} \left\{\inf_{\bar{\u} \in \mathcal{S}^{d-1},\bar{\u} \bot \u} P(\bar{\u}^\top X \leq \bar{\u}^\top \theta_0)\right\}\\
  & = & \inf_{\u\in \mathcal{S}^{d - 1}} P(\u^\top X \leq \u^\top \theta_0)\\
  & = & D(\theta_0, F).
\end{eqnarray*}
This proves $\varepsilon(T^*, \mathcal{X}^n) \stackrel{\text{a.s.}}{\longrightarrow} \frac{D(\theta_0, F)}{1 + D(\theta_0, F)} = 1/3$, because $D(\theta_0, F) = \frac{1}{2}$ under Conditions (C2)-(C3). This completes the proof of this theorem. \hfill $\Box$
\medskip

\noindent
\textbf{Remark 2}~~
It is worth noting that, 
both \emph{halfspace symmetry} and \emph{weak smooth condition} assumptions in this paper can \emph{not} be further relaxed if one wants to obtain \emph{exactly} the limiting breakdown point of \emph{HM}. The former is the weakest assumotion to guarantee to have a unique center. The latter is equivalent to the continuity of $D(\x, F)$ at $\theta_0$, which is necessary for deriving the limit for  both the lower and upper bound, while the upper bound given in Lemma 4 could not be further improved for fixed $n$. 

\vskip 0.1 in
\section{Concluding remarks}
\paragraph{}
\vskip 0.1 in \label{Conclusion}

In this paper, we consider the limit of the finite sample breakdown point of \emph{HM} under \emph{weaker} assumption on underlying distribution and data set.
Under such assumptions, the random observations may  \emph{not} be `in general position'. This causes additional inconvenience to the derivation of the limiting result compared to the scenario of $\mathcal{X}^n$ being in general position. During our investigation, relationships between various smooth conditions have been established and the representation of the Tukey depth and median regions has also been obtained without imposing the `in general position' assumption.
\medskip

Tukey halfspace depth idea has been extended beyond the location setting to many other settings (e.g., regression, functional data, etc.). We anticipate that our results here could also be extended to those settings.

\vskip 0.1 in
\section*{Acknowledgements}
\paragraph{}
\vskip 0.1 in \label{Conclusion}

The research of the first two authors is supported by National Natural Science Foundation of China (Grant No.11461029, 61263014, 61563018), NSF of Jiangxi Province (No.20142BAB211014, 20143ACB21012,  20132BAB201011, 20151BAB211016), and the Key Science Fund Project of Jiangxi provincial education department (No.GJJ150439, KJLD13033, KJLD14034).

\end{document}